\documentclass[a4paper,12pt]{article}

\usepackage{latexsym}
\usepackage{graphicx}
\usepackage{amsmath, amsthm}
\usepackage{hyperref}

\newcommand{\mbold}[1]{ \mbox{\boldmath $#1$} }

\newtheorem{remark}{Remark}
\newtheorem{example}{Example}
\newtheorem{definition}{Definition}

\newtheorem{theorem}{Theorem}

 \title{\textbf{Convergence Criteria for a Hopfield-type Neural Network}
 }
\author{
Raveen Goundar\thanks{{Corresponding Author. Email: raveen\underline{ }g@hotmail.com}}\\ and \\Jito Vanualailai\\
{\small \em Department of Mathematics and Computing Science,}
\\
{\small \em University of the South Pacific, Suva, Fiji.}
\\
\date{}
\\
{\textmd{2000 Mathematics Subject Classification: 34D20, 92B20.}}
}
\begin{document}
\maketitle
\begin{abstract}
Motivated by recent applications of the Lyapunov's method in
artificial neural networks, which could be considered as dynamical
systems for which the convergence of the system trajectories to
equilibrium states is a necessity. We re-look at a well-known
Krasovskii's stability criterion pertaining to a non linear
autonomous system. Instead, we consider the components of the same
autonomous system with the help of the elements of Jacobian matrix
{\bf{J}}({\bf{x}}), thus proposing much simpler convergence
criteria via the method of Lyapunov. We then apply our results to
artificial neural networks and discuss our results with respect to
recent ones in the field.\\
\\
{\em Keywords and Phrases}: Lyapunov Stability, Hopfield-Tank
Neural Networks\\
\end{abstract}
\section{Introduction}
The {\em Direct Method of Lyapunov}, which utilizes energy-like
functions  called {\em Lyapunov functions}, is now a
well-entrenched technique in the qualitative analysis of
mathematical systems governed by differential equations. A flurry
of activities by mathematicians,  particularly within the period
of early 1940s and the late 1960s, extended the work of Lyapunov
to produce results that are now indispensable in many
applications. (A good modern review of the Lyapunov method and its
many applications is by Sastry~\cite{sastry}.) This paper is
motivated to a large extent by modern applications of the Lyapunov
method, especially in the field of artificial neural networks.

We start by considering the autonomous system  of the form
\begin{equation}
\mbold{x}'(t) = \mbold{g}(\mbold{x}) \,, \mbold{x}(t_{0}) =
\mbold{x}_{0} \,.
\label{auto1}
\end{equation}
Throughout the paper, guided by a well-known result of Krasovskii,
we will strive to portray a simple and flexible method of
proposing a stability criterion for system~(\ref{auto1}). We
conclude by considering an application in artificial neural
networks.

 Throughout the article, we suppose that, in system~(\ref{auto1}),
 $\mbold{g} \in C[ \mbold{R}^{n}, \mbold{R}^{n}]$, and
 is
 smooth enough to guarantee existence, uniqueness and continuous dependence
 of solutions $\mbold{x}(t)=\mbold{x}(t;  \mbold{x}_{0})$, with
 $\mbold{x}=(x_{1}, \ldots, x_{n})^{T}$. The following definition and theorems
 of Lyapunov will be used in this article. (We will
 use those in Glendenning~\cite{Glendenning}).
\begin{definition}
Suppose that the origin, ${\mbold x} = {\mbold 0}$, is an
equilibrium point for system~(\ref{auto1}). Let $D$ be an open
neighborhood of ${\mbold 0}$ and $ V : D \rightarrow {\mbold
 R}$ be a continuously differentiable function. Then we can define
 the derivative of $ V $ along trajectories by differentiating $V$
 with respect to time using the chain rule, so
 \[
V'({\mbold x}) = \frac{dV({\mbold x})}{dt} = {\mbold
x}'\cdot\nabla V({\mbold x}) = {\mbold g}({\mbold x})\cdot\nabla
V({\mbold x}) = \sum_{i=1}^{n} g_{i}({\mbold x}) \frac{\partial
V({\mbold x})}{\partial x_{i}}\,,
\]
where the subscripts denote the components of ${\mbold g}$ and
${\mbold x}$. Then $ V $ is a Lyapunov function on $ D $ iff
\begin{description}
\item[{\rm (i)}]$V$ is continuously differentiable on $D$\,;
\item[{\rm (ii)}] $V({\mbold 0}) = 0$ and $V({\mbold x}) > 0$ for all
${\mbold x}\in D \setminus \{ {\mbold 0}\}$\,;
\item[{\rm (iii)}] $V'({\mbold x}) \leq 0$ for all ${\mbold x} \in D$.
\end{description}
\label{c1d1}
\end{definition}
\begin{theorem}[Lyapunov's Stability Theorem]
Let ${\mbold x} = {\mbold 0}$ be an equilibrium point for
system~(\ref{auto1}) and $D \subset {\mbold R}^{n}$ be a domain
containing ${\mbold x} = {\mbold 0}$. Let $V({\mbold x})$ be a
Lyapunov function on an open neighborhood of $D$, then ${\mbold x}
= {\mbold 0}$ is stable.
\label{c1t1}
\end{theorem}
\begin{theorem}[Lyapunov's Asymptotic Stability Theorem]
Let ${\mbold x} = {\mbold 0}$ be an equilibrium point for
system~(\ref{auto1}) and $D \subset {\mbold R}^{n}$ be a domain
containing ${\mbold x} = {\mbold 0}$. Let $V({\mbold x})$ be a
Lyapunov function on an open neighborhood of $D$. If $V'({\mbold
0}) = 0$ and $V'({\mbold x}) < {\mbold 0}$ for all ${\mbold x} \in
D \setminus \{{\mbold 0}\}$, then ${\mbold x} = {\mbold 0}$ is
asymptotically stable.
\label{c1t2}
\end{theorem}
\begin{theorem}[Lyapunov's Theorem of Global Asymptotic Stability]
Let ${\mbold x} = {\mbold 0}$ be a equilibrium point for
system~(\ref{auto1}) and let $V({\mbold x})$ be a Lyapunov
function for all ${\mbold x} \in {\mbold R}^{n}$. If ${\mbold x} =
{\mbold 0}$ is asymptotically stable and $V({\mbold x})$ is
radially unbounded, then ${\mbold x} = {\mbold 0}$ is globally
asymptotically stable.
\label{c1t3}
\end{theorem}
We carry the assumption that $\mbold{g}({\bf 0}) \equiv {\bf 0}$
so that ${\bf 0}$ is the zero solution of (\ref{auto1}).

\section{Convergence Criteria}
In 1954, Krasovskii~\cite{kra} established an asymptotic stability
criterion that avoided the linearization principle, and in the
process established a method of estimating the extent of
asymptotic stability region for a nonlinear systems. He assumed
that $\mbold{g} \in C^\prime[\mbold{R}^{n}, \mbold{R}^{n}]$ and
$\mbold{g}({\bf 0}) = {\bf 0}$. Then system~(\ref{auto1}) can be
written as
\[
\mbold{x}'(t) = \int_{0}^{1} \mbold{J}(s\mbold{x}) \mbold{x}  ds
\]
where $\mbold{J}$ is the Jacobian matrix
\[
\mbold{J}(\mbold{x}) = \frac{\partial
\mbold{g}(\mbold{x})}{\partial \mbold{x}} \,.
\]
The following result
by Krasovskii is a fundamental one in control theory.
\begin{theorem}[Krasovskii~\cite{kra}] Let $\mbold{g} \in C^{\prime}[\mbold{R}^{n}, \mbold{R}^{n}]$
and $\mbold{g}({\bf 0}) = {\bf 0}$. If there exists a constant
positive definite symmetric matrix ${\bf P}$ such that
\[
\mbold{x}^{T} [ \mbold{P} \mbold{J}(\mbold{x}) +
\mbold{J}^{T}(\mbold{x}) \mbold{P} ] \mbold{x}
\]
is a negative definite function, then the zero solution of
system~(\ref{auto1}) is globally asymptotically stable.
\end{theorem}
For our purpose, we need a criterion that explicitly uses each
component of system~(\ref{auto1}). Thus, using the elements of
Jacobian matrix; $J_{ij}({\mbold{x}})$, we define
\begin{equation}
\mbold{D}(\mbold{x}) = \left[ d_{ij}(\mbold{x}) \right]_{ n \times
n }
\label{c1e2}
\end{equation}
\\
where \[d_{ij}(\mbold{x}) =\int_{0}^{1}J_{ij}(s\mbold{x})ds =
\int_{0}^{1} \frac{ \partial g_{i}(s \mbold{x}) }{ \partial
(sx_{j}) } ds \,,\]\\
such that system~(\ref{auto1}) can be rewritten as
\begin{equation}
\mbold{x}'(t) = \mbold{D}(\mbold{x}) \mbold{x} \,.
\label{form}
\end{equation} A decoupled form for the $i$-th component of
system~(\ref{form}) is
\begin{equation}
x_{i}'(t) = d_{ii}(\mbold{x}) x_{i} + \sum_{{j=1} \atop  {j \neq
i}}^{n} d_{ij}(\mbold{x}) x_{j} \,.
\label{c2e2}
\end{equation}
\begin{remark}
{\em Note that in (\ref{c2e2}), the term $d_{ij}(\mbold{x})x_{j}$,
for $i, j = 1, \ldots, n$, is continuously differentiable with
respect to $\mbold{x} \in \mbold{R}^{n}$ for the simple reason
that $\mbold{D(x)x=g(x)}$ and $\mbold{g} \in
C^{\prime}[\mbold{R}^{n}, \mbold{R}^{n}]$. \label{condiff} }
\end{remark}
The following result of ours, guarantees the convergence criteria
for autonomous system~(\ref{auto1}).

\begin{theorem}
Let $\mbold{g} \in C^{\prime}[\mbold{R}^{n}, \mbold{R}^{n}]$ and
$\mbold{g}({\bf 0}) = {\bf 0}$. Let
\[
\beta_{i}(\mbold{x}) =
d_{ii}(\mbold{x}) + \frac{1}{2}
\sum_{{j=1}
\atop  {j \neq i}}^{n}
\left( | d_{ij}(\mbold{x}) |  + | d_{ji}(\mbold{x}) | \right) \,.
\]
Define $D = \{{\mbold x} \in {\mbold R}^{n} : \|{\mbold x}\| \leq
M \}$ for some $M > 0$ and assume that $d_{ij}(\mbold{x})x_{i}$
are continuous on $\mbold{R}^{n}$ for $i,j=1, \ldots,  n$, such
that $i \neq j$. Then the zero solution of (\ref{auto1}) is

\begin{enumerate}
\item [(a)] stable if $-\infty < \beta_{i}({\mbold x})
\leq 0$ for $i = 1,2, \ldots ,n$ and ${\mbold x} \in D$.
\item [(b)] asymptotically stable if $-\infty <
\beta_{i}({\mbold x}) < 0$ for $i = 1,2,\ldots,n$ and ${\mbold x}
\in D$.
\item [(c)] globally asymptotically stable if $-\infty <
\beta_{i}({\mbold x}) < 0$ for all ${\mbold x} \in {\mbold
R}^{n}$.
\end{enumerate}
\label{raveen}
\end{theorem}

\begin{proof}
Consider
\[
V(\mbold{x}) = \frac{1}{2} \sum_{i=1}^{n} x_{i}^{2}
\]
as a tentative Lyapunov function  for system~(\ref{auto1}). We
have, along a solution of (\ref{auto1}),
\begin{eqnarray}
\frac{d}{dt} \left[ V \right]_{(\ref{auto1})}
& =  & \frac{1}{2} \sum_{i=1}^{n} \frac{d}{dt}
\left[ x_{i}^{2} \right]
 =
\sum_{i=1}^{n}  x_{i} x_{i}'(t)\nonumber \\
& = &  \sum_{i=1}^{n}   x_{i} \left[
d_{ii}(\mbold{x})x_{i}  +
\sum_{{j=1} \atop  {j \neq i}}^{n} d_{ij}(\mbold{x}) x_{j}  \right]\nonumber \\
& = & \sum_{i=1}^{n}   \left[
d_{ii}(\mbold{x}) x_{i}^{2}
+   \sum_{{j=1}
\atop  {j \neq i}}^{n} d_{ij}(\mbold{x}) x_{j} x_{i} \right]\nonumber \\
& = & \sum_{i=1}^{n}   \left[
d_{ii}(\mbold{x}) x_{i}^{2}
+   \frac{1}{2} \sum_{{j=1}
\atop  {j \neq i}}^{n} [ d_{ij}(\mbold{x}) + d_{ji}(\mbold{x}) ] x_{j} x_{i} \right] \nonumber\\
& \leq & \sum_{i=1}^{n}   \left[
d_{ii}(\mbold{x}) x_{i}^{2}
+   \frac{1}{4} \sum_{{j=1}
\atop  {j \neq i}}^{n} [ | d_{ij}(\mbold{x})| + |d_{ji}(\mbold{x})| ]
(x_{j}^{2} + x_{i}^{2}) \right]\nonumber \\
& = & \sum_{i=1}^{n}   \left[
d_{ii}(\mbold{x})
+
\frac{1}{2} \sum_{{j=1}
\atop  {j \neq i}}^{n}
[ | d_{ij}(\mbold{x})| + |d_{ji}(\mbold{x})| ]
 \right] x_{i}^{2} \label{rav1}\\
& = & \sum_{i=1}^{n} \beta_{i}(\mbold{x}) x_{i}^{2} \,.
\label{c1e3}
\end{eqnarray}
Expanded form of system~(\ref{rav1}) is
\[
\frac{dV}{dt} \leq \sum_{i=1}^{n} \left[ d_{ii}(\mbold{x}) x_{i} x_{i}
+
\frac{1}{2} \sum_{{j=1}
\atop  {j \neq i}}^{n}
[ | d_{ij}(\mbold{x})x_{i}x_{i}| + |d_{ji}(\mbold{x})x_{i}x_{i}| ]
\right] \,.
\]
By Remark~\ref{condiff}, the first and third terms of
system~(\ref{rav1}) are continuous on $\mbold{R}^{n}$, and by
assumption of Theorem~\ref{raveen}, the second term is also
continuous on $\mbold{R}^{n}$. Hence $V({\mbold x})$ is continuous
on ${\mbold R}^{n}$. Since
\begin{equation}
V({\mbold x}) = \frac{1}{2} \sum_{i=1}^{n}{x_{i}}^{2}\,,
\label{c1e4}
\end{equation}
we have therefore, $V({\mbold 0}) = 0$ and $V({\mbold x}) > 0$ for
all ${\mbold x} \in {\mbold R}^{n}\setminus\{\mbold{0}\}$. From
equation~(\ref{c1e3}),
\begin{equation}
V'({\mbold x}) \leq \sum_{i=1}^{n} \beta_{i}({\mbold
x}){x_{i}}^{2} \label{c1e5}
\end{equation}
and by condition~(a) of Theorem~\ref{raveen}, we have $V'({\mbold
x}) \leq 0$ for all ${\mbold x} \in D$. Hence by
Theorem~\ref{c1t1}, the zero solution of system~(\ref{auto1}) is
stable. Moreover, by condition~(b) of Theorem~\ref{raveen},
equation~(\ref{c1e5}) implies $V'({\mbold 0}) = 0$ and $V'({\mbold
x}) < 0$ for all ${\mbold x} \in D \setminus \{ {\mbold 0}\}$.
Hence by Theorem~\ref{c1t2}, the zero solution of
system~(\ref{auto1}) is asymptotically stable. Furthermore, by
condition~(c) of Theorem~\ref{raveen}, equation~(\ref{c1e5})
implies $V'({\mbold 0}) = 0$ and $V'({\mbold x}) < 0$ for all
${\mbold x} \in {\mbold R}^{n}$. Note that (\ref{c1e4}) implies
$V({\mbold x}) \rightarrow \infty$ as $\|{\mbold x}\| \rightarrow
\infty$, thus $V({\mbold x})$ is radially unbounded. Hence by
Theorem~\ref{c1t3}, the zero solution of system~(\ref{auto1}) is
globally asymptotically stable.
\end{proof}
Let us consider some examples to show the applicability of
Theorem~\ref{raveen}.
\begin{example}
{\em We consider the following two-dimensional system
\begin{eqnarray}
\left[
\begin{array}{c}
x_{1}'(t)\\
x_{2}'(t)
\end{array}
\right]
 =
\left[
\begin{array}{c}
-2 x_{1} + x_{2}^{2}\\
x_{1}^{2} - 2 x_{2}
\end{array}
\right]\,,
\label{c2e8}
\end{eqnarray}
with $x_{1}(t_{0}) = x_{10}$ and $x_{2}(t_{0}) = x_{20}$. In the
form of system~(\ref{form}), system~(\ref{c2e8}) can be written as
\begin{eqnarray*}
\left[
\begin{array}{c}
x_{1}'(t)\\
x_{2}'(t)
\end{array}
\right]  = \left[
\begin{array}{cc}
-2 & x_{2}\\
x_{1} & -2
\end{array}
\right] \left[
\begin{array}{c}
x_{1}\\
x_{2}
\end{array}
\right]\,.
\end{eqnarray*}
}
\label{c2ex1}
\end{example}
The assumption of Theorem~\ref{raveen} is satisfied since
\begin{equation}
d_{12}({\mbold x})x_{1}  = d_{21}({\mbold x})x_{2} = x_{1}x_{2}\,.
\nonumber
\end{equation}
Next we shall check condition~(a) of Theorem~\ref{raveen}. We have
\begin{eqnarray*}
\beta_{1}({\mbold x}) & = & d_{11}({\mbold x}) +
\frac{1}{2}\left(|d_{12}({\mbold x})| + |d_{21}({\mbold
x})|\right)
\\
& = & -2 + \frac{1}{2}\left(|x_{2}| + |x_{1}|\right)\nonumber\,.
\end{eqnarray*}
Solving the inequality $\beta _{1}({\mbold x}) < 0$, we have
\[
|x_{1}| + |x_{2}| < 4\,,
\]
and `squaring' both sides gives
\[
{x_{1}}^{2} + {x_{2}}^{2} + 2|x_{1}||x_{2}| < 16\,.
\]
Now
\[
{x_{1}}^{2} + {x_{2}}^{2} + 2|x_{1}||x_{2}| < {x_{1}}^{2} +
{x_{2}}^{2} + 2 \times \frac{1}{2}\left({x_{1}}^{2} +
{x_{2}}^{2}\right) = 2{x_{1}}^{2} + 2{x_{2}}^{2}\,.
\]
Then let
\[
2{x_{1}}^{2} + 2{x_{2}}^{2} < 16
\]
from which
\begin{equation}
{x_{1}}^{2} + {x_{2}}^{2} < 8\,.
\nonumber
\end{equation}
Similarly solving $\beta _{2}({\mbold x}) < 0$, we have
\begin{equation}
\beta_{2}({\mbold x}) = d_{22}({\mbold x}) +
\frac{1}{2}\left(|d_{21}({\mbold x})| + |d_{12}({\mbold
x})|\right) < 0\,,
\nonumber
\end{equation}
which gives
\begin{equation}
-2 + \frac{1}{2} \left(|x_{1}| + |x_{2}|\right) < 0.
\label{rav2}
\end{equation}
Further simplification of (\ref{rav2}) gives us
\begin{equation}
{x_{1}}^{2} + {x_{2}}^{2} < 8\,.
\nonumber
\end{equation}
Therefore, let
\begin{equation}
D = \{{\mbold x} \in {\mbold R}^{2} : \|{\mbold x}\| <
\sqrt{8}\}\,.
\nonumber
\end{equation}
Hence by condition~(a) of Theorem~\ref{raveen}, the zero solution
of system~(\ref{c2e8}) is asymptotically stable.
\begin{example}
{\em We consider the following two-dimensional system
\begin{eqnarray}
\left[
\begin{array}{c}
x_{1}'(t)\\
\\
x_{2}'(t)
\end{array}
\right] = \left[
\begin{array}{c}
\displaystyle -4 x_{1} +  x_{1}\mbox{sech}(x_{1}) + 4 x_{2}\\
\\\displaystyle -x_{1} - 6 x_{2} - x_{2}\cos{(x_{2})}
\end{array}
\right]\,,
\label{c2e12}
\end{eqnarray}
which can be written in the form of system~(\ref{form}) as
\begin{eqnarray*}
\left[
\begin{array}{c}
x_{1}'(t)\\
\\
x_{2}'(t)
\end{array}
\right] = \left[
\begin{array}{cc}
\displaystyle -4 + \mbox{sech}(x_{1}) & \displaystyle 4\\
\\
\displaystyle -1 & \displaystyle -6 - \cos(x_{2})
\end{array}
\right] \left[
\begin{array}{c}
x_{1}\\
\\
x_{2}
\end{array}
\right]\,.
\end{eqnarray*}
The assumption of Theorem~\ref{raveen} is satisfied since
$d_{12}({\mbold x}) x_{1} = 4x_{1}$ and $d_{21}({\mbold x})x_{2} =
-x_{2}$. Next we shall check condition~(c) of
Theorem~\ref{raveen}. We have
\begin{eqnarray*}
\beta_{1}({\mbold x}) & = & d_{11}({\mbold x}) + \frac{1}{2}\left(
|d_{12}({\mbold x})| + |d_{21}({\mbold x})|\right)\\
& = &
 -4 + \mbox{sech}(x_{1}) + \frac{1}{2}\left(|4| +
|-1|\right)\\
& = & \mbox{sech}(x_{1}) - \frac{3}{2} \leq 1 - \frac{3}{2} =
-\frac{1}{2} < 0\,.
\end{eqnarray*}
Similarly, we have
\begin{eqnarray*}
\beta_{2}({\mbold x}) & = & d_{22}({\mbold x}) +
\frac{1}{2}\left(|d_{21}({\mbold
x})| + |d_{12}({\mbold x})|\right)\\
& = & -6  - \cos(x_{2}) + \frac{1}{2}\left(|-1| + |4|\right)\\
& = & -\frac{7}{2} - \cos(x_{2}) \leq -\frac{7}{2} + 1 =
-\frac{5}{2} < 0\,.
\end{eqnarray*}
Since both $\beta_{1}({\mbold x}) < 0$ and $\beta_{2}({\mbold x})
< 0$ for all ${\mbold x} \in {\mbold R}^{2}$ hence by
condition~(c) of Theorem~\ref{raveen}, the zero solution of
system~(\ref{c2e12}) is globally asymptotically stable. }
\end{example}
\section{Application in Artificial Neural Networks}
Artificial neural networks (ANNs) can be considered as dynamical
systems with several equilibrium states. An essential operating
condition for a neural network is that all system trajectories
must converge to the equilibrium states. (A good overview of the
concepts associated with biological neural networks is given in
\cite{arbib}).

We will consider an ANN that is described thoroughly in
Lakshmikantham~et~al.~\cite{laks1}, and provide a stability
criteria using Theorem~\ref{raveen}. The ANN in question has $n$
units. To the $i$th unit, we associate its {\em activation
state\/} at time $t$, a real number $x_{i}=x_{i}(t)$; an {\em
output function\/} $\mu_{i}$; a fixed {\em bias\/} $\theta_{i}$;
and an {\em output signal\/} $R_{i}=\mu_{i}(x_{i}+\theta_{i})$.
The {\em weight\/} or connection strength on the line from unit
$j$ to unit $i$ is a fixed real number $W_{ij}$. When $W_{ij}=0$,
there is no transmission from unit $j$ to unit $i$. The {\em
incoming signal\/} from unit $j$ to unit $i$ is
$S_{ij}=W_{ij}R_{j}$. In addition, there can be a vector
$\mbold{I}$ of any number of {\em external inputs\/} feeding into
some or all units, so that we may write $\mbold{I}=(I_{1}, \ldots,
I_{m})^{T}$.

An ANN with fixed weights is a dynamical system: given initial values
of the activation of all units, the future activations can be computed. The
future activation states are assumed to be determined by a system of $n$
differential equations, the $i$th equation of which is
\begin{eqnarray}
x_{i}'(t) & = &
G_{i}(x_{i}, S_{i1}, \ldots, S_{in}, \mbold{I})  =
G_{i}(x_{i}, W_{i1}R_{1}, \ldots, W_{in}R_{n}, \mbold{I})
\nonumber \\
& = &
G_{i}(x_{i}; W_{i1} \mu_{1}(x_{1}+\theta_{1}), \ldots,
W_{in} \mu_{n}(x_{n}+\theta_{n});
I_{1}, \ldots, I_{m} ) \,.
\label{net1}
\end{eqnarray}
With $W_{ij}$, $\theta_{i}$ and $I_{k}$ assumed known, we can write
(\ref{net1}) as
\begin{equation}
x_{i}'(t)= g_{i}(x_{1}, \ldots, x_{n}) \,,
\label{net2}
\end{equation}
or in vector notation
\begin{equation}
\mbold{x}'(t) = \mbold{g}(\mbold{x}) \,,
\label{net3}
\end{equation}
where $\mbold{g}$ is a vector on Euclidean space $\mbold{R}^{n}$
whose $i$th element is $g_{i}$ given in (\ref{net2}). We assume
that $\mbold{g}$ is continuously differentiable and satisfies the
usual theorems on existence, continuity and uniqueness of
solutions. Thus, since $\mbold{g} \in C^{\prime}[\mbold{R}^{n},
\mbold{R}^{n}]$, we can define $\mbold{D}(\mbold{x})$ as in
(\ref{c1e2}) but using $\mbold{g}$ in (\ref{net3}). Hence, if
$\mbold{g}(\bf 0) \equiv {\bf 0}$, then system~(\ref{net3}) can be
written as
\[
\mbold{x}'(t) = \mbold{D}(\mbold{x}) \mbold{x} \,, \;\;
\mbold{x}(t_{0}) = \mbold{x}_{0} \,,
\]
the $i$th component of which in a decoupled form is
\[
x_{i}'(t) = d_{ii}(\mbold{x}) x_{i}
+
\sum_{{j=1}
\atop  {j \neq i}}^{n}
d_{ij}(\mbold{x}) x_{j} \,.
\]

First, we state a comparable result by Lakshmikantham~et~al.~\cite{laks1},
page~152, who used the
concept of {\em vector Lyapunov functions}.
\begin{theorem}[Lakshmikantham, Matrosov and Sivasundaram~\cite{laks1}]
Let $\mbold{g} \in C^{\prime}[\mbold{R}^{n}, \mbold{R}^{n}]$ and
$\mbold{g}({\bf 0}) = {\bf 0}$. Let
\begin{equation}
\beta_{i}(\mbold{x}) =
d_{ii}(\mbold{x}) +
\sum_{{j=1}
\atop  {j \neq i}}^{n}
| d_{ij}(\mbold{x}) |  \,.
\label{laks1}
\end{equation}
Suppose that
\begin{equation}
\beta_{i}(\mbold{x})  < 0 \;\; \mbox{ if }
\;\; x_{i}^{2} \geq x_{j}^{2}  \,,
\label{laks}
\end{equation}
for
$i, j =1, \ldots, n$ and $\mbold{x} \in \mbold{R}^{n}$,
$\mbold{x} \neq {\bf 0}$. Then the
zero solution of (\ref{net3}) is globally asymptotically
stable.
\label{LaksThm1}
\end{theorem}
If we apply condition~(b) of Theorem~\ref{raveen}, then we obtain
a simpler convergence criteria.
\begin{theorem}
Let $\mbold{g} \in C^{\prime}[\mbold{R}^{n}, \mbold{R}^{n}]$ and
$\mbold{g}({\bf 0}) = {\bf 0}$. Let
\[
\beta_{i}(\mbold{x}) = d_{ii}(\mbold{x}) + \frac{1}{2} \sum_{{j=1}
\atop  {j \neq i}}^{n} \left( | d_{ij}(\mbold{x}) |  + |
d_{ji}(\mbold{x}) | \right) \,.
\]
Define $D = \{{\mbold x} \in {\mbold R}^{n} : \|{\mbold x}\| \leq
M \}$ for some $M > 0$ and assume that $d_{ij}(\mbold{x})x_{i}$
are continuous on $\mbold{R}^{n}$ for $i,j=1, \ldots,  n$, such
that $i \neq j$. Then the zero solution of (\ref{net3}) is
asymptotically stable if $-\infty < \beta_{i}({\mbold x}) < 0$ for
$i = 1,2,\ldots,n$ and ${\mbold x} \in D$. \label{LaksThm11}
\end{theorem}
Thus, the application of Theorem~\ref{raveen} to artificial neural
network, considering system~(\ref{net3}), gives us a simpler
criterion guaranteeing asymptotic stability as showed by
Theorem~\ref{LaksThm1}. Hence the strong condition $ x_{i}^{2}
\geq x_{j}^{2}$ that appears in Theorem~\ref{LaksThm1} is not
necessary.

Next, we look at a specific case of (\ref{net3}). The specific ANN
is of the additive type and is often referred to as the
Hopfield-Tank ANN, a much studied class of network
dynamics~\cite{hopfield}. It is described by the nonlinear
differential equation
\begin{eqnarray}
x_{i}'(t)
& = & -a_{i} x_{i}(t) +  \sum_{j=1}^{n} W_{ij} \; \mu_{j}(x_{j}(t)
+ \theta_{j}) +
I_{i}(t) \nonumber \\
& = & -a_{i} x_{i}(t) +  \sum_{j=1}^{n} W_{ij} \; \nu_{j}(x_{j}(t)) +
I_{i}(t) \,,
\label{NN}
\end{eqnarray}
where  $a_{i}>0$ is the constant {\em decay rate\/},
$I_{i}(t)$ is the external input (to the $i$th neuron)
defined almost everywhere on $[0, \infty)$ and $\nu_{i}$ is the
suppressed notation for the
fixed $\theta_{i}$ by having $\theta_{i}$ incorporated into
$\mu_{i}$. The function $\nu_{i}$ is called the
{\em neuron activation function\/}.

Now, define $\mbold{A}=\mbox{diag}(-a_{1}, \ldots, -a_{n})$,
$\mbold{x}=
(x_{1}, \ldots , x_{n})^{T}$,
\[
h_{i}(\mbold{x})= \sum_{j=1}^{n}  W_{ij} \nu_{j}(x_{j}) \mbox{
with } \mbold{h}(\mbold{x})=(h_{1}(\mbold{x}), \ldots,
h_{n}(\mbold{x}))^{T} \,,
\]
and $\mbold{u}(t) = (I_{i}(t), \ldots, I_{n}(t))^{T}$. Then
(\ref{NN}) is the $i$th component of the system
\begin{equation}
\mbold{x}'(t) = \mbold{A} \mbold{x} + \mbold{h}(\mbold{x}) +
\mbold{u}(t) \,, \;\; \mbold{x}(t_{0}) = \mbold{x}_{0} \,.
\label{input}
\end{equation}
When the external input vector, $\mbold{u}$, is zero, the
nonautonomous system (\ref{input}) reduces to the autonomous
system
\begin{equation}
\mbold{x}'(t) = \mbold{A} \mbold{x} + \mbold{h}(\mbold{x}) \,,
\;\; \mbold{x}(t_{0}) = \mbold{x}_{0} \,. \label{noinput}
\end{equation}
For this, we assume that $\mbold{x}^{*}=(x_{1}^{*}, \ldots,
x_{n}^{*})^{T}$ is an equilibrium point, so that $\mbold{A}
\mbold{x}^{*} + \mbold{h}(\mbold{x}^{*}) = {\bf 0}$. By
translating the origin, ${\bf 0}$, to this equilibrium point, we
can make ${\bf 0}$ an equilibrium point. In this case,
$\mbold{h}({\bf 0}) \equiv {\bf 0}$.  Since this is of great
notational help, we will henceforth consider ${\bf 0}$ as an
equilibrium point or zero solution of (\ref{noinput}).

Let us next assumed that $\mbold{h} \in C^{\prime}[\mbold{R}^{n},
\mbold{R}^{n}]$. Then using the elements of Jacobian matrix,
$J_{ij}(\mbold{x})$, we define
\[
\mbold{F}(\mbold{x}) = \left[ f_{ij}(\mbold{x}) \right]_{ n \times
n }  \;\;  \mbox{ where }  \;\; f_{ij}(\mbold{x}) =
\int_{0}^{1}J_{ij}(\mbold{x})ds = \int_{0}^{1} \frac{
\partial h_{i}(s \mbold{x}) }{
\partial (sx_{j}) } ds ,
\]
hence system~(\ref{noinput}) can be rewritten as
\begin{equation}
\mbold{x}'(t) = \mbold{A}\mbold{x} + \mbold{F}(\mbold{x})
\mbold{x} = \left[ \mbold{A} + \mbold{F}(\mbold{x}) \right]
\mbold{x} \,.
\label{noinput2}
\end{equation}
The $i$th component of (\ref{noinput2}) in a decoupled form is
\[
x_{i}'(t) = [ -a_{i} + f_{ii}(\mbold{x})] x_{i}(t) + \sum_{{j=1}
\atop {j \neq i}}^{n}  f_{ij}(\mbold{x}) x_{j} \,.
\]
Thus the following theorem is an application of our result;
Theorem~\ref{raveen}.
\begin{theorem}
Let $\mbold{h} \in C^{\prime}[\mbold{R}^{n}, \mbold{R}^{n}]$ and
$\mbold{h}({\bf 0}) = {\bf 0}$. Let
\[
\beta_{i}(\mbold{x}) =  -a_{ii} + f_{ii}(\mbold{x}) + \frac{1}{2}
\sum_{{j=1} \atop  {j \neq i}}^{n} \left( | f_{ij}(\mbold{x}) |  +
| f_{ji}(\mbold{x}) | \right) \,.
\]
Define $D = \{{\mbold x} \in {\mbold R}^{n} : \|{\mbold x}\| \leq
M \}$ for some $M > 0$ and assume that $f_{ij}(\mbold{x})x_{i}$
are continuous on $\mbold{R}^{n}$ for $i,j=1, \ldots,  n$, such
that $i \neq j$. Then the zero solution of (\ref{noinput}) is

\begin{enumerate}
\item [(a)] stable if $-\infty < \beta_{i}({\mbold x})
\leq 0$ for $i = 1,2, \ldots ,n$ and ${\mbold x} \in D$.
\item [(b)] asymptotically stable if $-\infty <
\beta_{i}({\mbold x}) < 0$ for $i = 1,2,\ldots,n$ and ${\mbold x}
\in D$.
\item [(c)] globally asymptotically stable if $-\infty <
\beta_{i}({\mbold x}) < 0$ for all ${\mbold x} \in {\mbold
R}^{n}$.
\end{enumerate}
\label{raveen-2}
\end{theorem}
\begin{proof}
Applying Theorem~\ref{raveen} to system~(\ref{noinput}), and hence
to system~(\ref{noinput2}), with $\mbold{D}(\mbold{x}) = \mbold{A}
+ \mbold{F}(\mbold{x})$, $d_{ii}(\mbold{x})=-a_{i} +
f_{ii}(\mbold{x})$ and $d_{ij}(\mbold{x}) = f_{ij}(\mbold{x})$, we
easily obtain the conclusion of Theorem~\ref{raveen-2}.
\end{proof}
Let us consider one example of Theorem~\ref{raveen-2}.
\begin{example}{\em
Let us consider two-neural autonomous system.
\begin{eqnarray}
\left[
\begin{array}{c}
x_{1}'(t)\\
x_{2}'(t)
\end{array}
\right] = \left[
\begin{array}{cc}
-a_{1} & 0\\
0 & -a_{2}
\end{array}
\right] \left[
\begin{array}{c}
x_{1}\\
x_{2}
\end{array}
\right] + \left[
\begin{array}{c}
h_{1}({\mbold x})\\
h_{2}({\mbold x})
\end{array}\right]
\label{c4e6}
\end{eqnarray}
}
\end{example}
with  $x_{1}(t_{0}) = x_{10},\,\,\,x_{2}(t_{0}) = x_{20},\,\,\,0
\leq t_{0} \leq t$,  where,
\begin{eqnarray*}
a_{1} & = & 10,\,\,\,a_{2} = 10\,,\\
h_{1}({\mbold x}) & = & B_{11}\nu_{1}(x_{1}) +
B_{12}\nu_{2}(x_{2}) =-3x_{1} + x_{2} - \tanh(3x_{1})\,,\\
h_{2}({\mbold x}) & = & B_{21}\nu_{1}(x_{1}) +
B_{22}\nu_{2}(x_{2}) = x_{1} - x_{2} + \frac{1}{5}\tanh(3x_{2})\,.
\end{eqnarray*}
In the form of system~(\ref{noinput2}), system~(\ref{c4e6}) can be
written as
\begin{eqnarray*}
\left[\begin{array}{c} x_{1}'(t)\\
x_{2}'(t)
\end{array}
\right] = \left(\left[
\begin{array}{cc}
-10 & 0 \\
0 & -10
\end{array}
\right] +\ \left[
\begin{array}{cc}
-3 - \tau(x_{1}(t)) & 1\\
1 & -1 + \displaystyle\frac{1}{5}\tau(x_{2}(t))
\end{array}\right]\right)\left[
\begin{array}{c}
x_{1}\\
x_{2}
\end{array}
\right]\,,
\end{eqnarray*}
where for $i = 1,2,$ we define
\[ \tau(x_{i}(t)) = \left\{
\begin{array}{ll}
\displaystyle \frac{\tanh(3x_{i})}{x_{i}} & \mbox{$x_{i} \neq 0$}\,,\\
3 & \mbox{$x_{i} = 0$}\,,
\end{array}
\right.
\]
noting that $ 0 < \tau(x_{i}) \leq 3$ for all $x_{i} \in {\mbold
R}^{2}$\,. The assumption of Theorem~\ref{raveen-2} is satisfied
since $f_{12}({\mbold x})x_{1} = x_{1}$ and $ f_{21}({\mbold
x})x_{2} = x_{2}$. Now we shall check condition~(c) of
Theorem~\ref{raveen-2}. We have
\begin{eqnarray}
\beta_{1}({\mbold x}) & = & -a_{1} + f_{11}({\mbold x}) +
\frac{1}{2}(|f_{12}({\mbold x})| + |f_{21}({\mbold x})|)\nonumber\\
& = &
-10 - 3 - \tau(x_{1}(t)) + \frac{1}{2}(|1| + |1|)\nonumber\\
& = & -12 - \tau(x_{1}(t))
\label{rat1}\\
& < & -12 \nonumber
\end{eqnarray}
for all ${\mbold x}\in{\mbold R}^{2}\setminus\{{\mbold 0}\}$ and
\begin{eqnarray}
\beta_{2}({\mbold x}) & = & -a_{2} + f_{22}({\mbold x}) +
\frac{1}{2}(|f_{21}({\mbold x})| + |f_{12}({\mbold x})|)\nonumber\\
& = & -10 - 1 + \frac{1}{5}\tau(x_{2}(t)) + \frac{1}{2}(|1| +
|1|)\nonumber\\
& = & -10 + \frac{1}{5}\tau(x_{2}(t))
\label{rat2}\\
& < & -10 + \frac{3}{5} = -\frac{47}{5}  \nonumber
\end{eqnarray}
for all ${\mbold x} \in {\mbold R}^{2}\setminus\{{\mbold 0}\}$.
Clearly, both $\beta_{1}({\mbold x}) < 0$ and $\beta_{2}({\mbold
x}) < 0$ for all ${\mbold x}\in{\mbold R}^{2}\setminus\{{\mbold
0}\}$.

Next, we shall check the condition on $\beta_{i}({\mbold x})$ for
${\mbold x} = {\mbold 0}$, where $ i = 1,2$. From (\ref{rat1}), we
have
\[
\beta_{1}({\mbold x})  =  -12 - \tau(x_{1}(t))\,.
\]
Therefore,
\[
\beta_{1}({\mbold 0})  = -12 - 3 = -15 \,.
\]
Similarly, from (\ref{rat2}), we have
\[
\beta_{2}({\mbold x}) = -10 + \frac{1}{5} \tau(x_{2}(t))\,.
\]
Therefore,
\[
\beta_{2}({\mbold 0}) = -10 + \frac{3}{5} = -\frac{47}{5} \,.
\]
Since $\beta_{1}({\mbold x}) < 0$ and $\beta_{2}({\mbold x}) < 0$
for all ${\mbold x} \in {\mbold R}^{2}$, therefore, by
condition~(c) of Theorem~\ref{raveen-2}, the zero solution of
system~(\ref{c4e6}) is globally asymptotically stable.
\section{Conclusion}
We have established the criteria for stability, asymptotic
stability and global asymptotic stability for a non linear
autonomous system via the method of Lyapunov. We have also
considered the usefulness of our main results by application of it
to artificial neural networks.

Further research in this direction is being carried out,
considering a non autonomous system, wherein the external input
source is not assumed to be zero. Determining the convergence
criteria for a non autonomous system and to measure its rate of
convergence will be of grandness in applications to artificial
neural networks.

\end{document}